\documentclass{amsart}

\usepackage{amssymb}
\usepackage[all]{xy}
\usepackage{hyperref}

\usepackage{enumitem}   

\setlist[enumerate]{itemsep=.2em,topsep=.2em,leftmargin=1.25em,itemindent=2.0em}

\newtheorem{thm}{Theorem}

\newtheorem{lem}[thm]{Lemma}
\newtheorem{cor}[thm]{Corollary}

\theoremstyle{definition}
\newtheorem{defn}[thm]{Definition}

\newtheorem{say}[thm]{}
\newtheorem{exmp}[thm]{Example}

\newtheorem{rem}[thm]{Remark}

\newtheorem*{ack}{Acknowledgments}      
\newtheorem{notation}[thm]{Notation}   
  
\newtheorem{defn-thm}[thm]{Definition--Theorem}  
\newtheorem{defn-lem}[thm]{Definition--Lemma}  

\theoremstyle{remark}

\let \cedilla =\c
\renewcommand{\c}[0]{{\mathbb C}}  

\renewcommand{\o}[0]{{\mathcal O}} 
\newcommand{\z}[0]{{\mathbb Z}}

\newcommand{\p}[0]{{\mathbb P}}

\newcommand{\q}[0]{{\mathbb Q}}
\newcommand{\map}[0]{\dasharrow}
\newcommand{\qtq}[1]{\quad\mbox{#1}\quad}

\newcommand{\pics}[0]{\operatorname{\mathbf{Pic}}}

\newcommand{\gal}[0]{\operatorname{Gal}}

\newcommand{\supp}[0]{\operatorname{Supp}}

\newcommand{\bs}[0]{\operatorname{Bs}}
\newcommand{\aut}[0]{\operatorname{Aut}}

\newcommand{\sing}[0]{\operatorname{Sing}}

\newcommand{\jac}[0]{\operatorname{\mathbf{Jac}}} 

\newcommand{\chr}[0]{\operatorname{char}}

\newcommand{\grass}[0]{\operatorname{Grass}}

\newcommand{\tsum}[0]{\textstyle{\sum}}

\def\into{\DOTSB\lhook\joinrel\to}

\def\loccoh#1.#2.#3.#4.{H^{#1}_{#2}(#3,#4)}

\DeclareMathAlphabet{\mathchanc}{OT1}{pzc}%
                                {m}{it}

\newcommand{\gm}[0]{{\mathbb G}_m}

\newcommand{\PGL}{\mathrm{PGL}}

\usepackage[all]{xy}\xyoption{dvips}

\newcommand{\tprod}[0]{\textstyle{\prod}}

\begin{document}

\bibliographystyle{amsalpha}


\title{Theorems of Bertini and Chevalley}
 \author{J\'anos Koll\'ar}

 \begin{abstract}
   We give a short proof of   Chevalley's theorem that every  algebraic group is an extension of  an Abelian variety by a linear algebraic group. Along the way we treat   Bertini's irreducibility theorem.

 \end{abstract}

 \maketitle

 We  revisit the following two  well-known theorems.

 \begin{thm}[Bertini's irreducibility theorem]\label{bert.irr.thm}
   Let $X$ be a geometrically   irreducible variety over an infinite field $k$, and
   $\pi:X\map \p^n$ a dominant morphism, $n\geq 2$.
   Let $H\subset \p^n$ be  a general hyperplane. Then
   $\pi^{-1}(H)$ is geometrically  irreducible.
   \end{thm}

 \begin{thm}[Chevalley or  Barsotti-Chevalley theorem]
   \label{chev.thm} Let $k$ be a perfect field and 
  $G$ a connected, smooth, algebraic group over $k$.
  Then $G$ can be written as  a group extension
  $$
  1\to (\mbox{affine algebraic group}) \to G\to
  (\mbox{proper algebraic group}) \to 0.
  $$
\end{thm}

 The aim is to give proofs that use only results from    Shafarevich's textbook \cite{shaf}.
The proof of Bertini's irreducibility theorem 
given in Section~\ref{bert.sec}  is such, and  seems shorter than other proofs in the literature.

For Chevalley's theorem we also need the
existence of the Jacobian of a smooth, projective curve; see Theorem~\ref{jac.exists}.
From algebraic group theory we use only Weil's theorems on rational maps to algebraic groups, these are reproved in Section~\ref{sec.p1.pf}.
See  Remark~\ref{chev.thm.rem.1} for other versions of Chevalley's theorem.

The main message is that one should not be afraid of working with Weil divisors; they provide the optimal framework for some questions.

\begin{say}[Historical comments and references]

  The  original statements of Bertini's  theorems and various proofs are discussed in \cite{MR1661859}.  While Bertini's smoothness  theorem is treated in many places, the irreducibility theorem
  appears less frequently in textbooks.
  (The `First Bertini theorem' in \cite[Sec.II.6.1]{shaf} is a  different result about irreducibility.)
  
  \cite[Sec.X.13]{hodge-ped}, \cite[Sec.I.6]{MR725671} and
   \cite[\href{https://stacks.math.columbia.edu/tag/0G4C}{Tag 0G4C}]{stacks-project} give  related  proofs of
   Theorem~\ref{bert.irr.thm}.
\cite[Sec.3.4]{Flenner-OCarrol-Vogel} and 
\cite[\href{https://stacks.math.columbia.edu/tag/0G4G}{Tag 0G4G}]{stacks-project}
   derive it from  Bertini-type results for normal
varieties (first fully treated  in \cite{MR37548}) and the
Enriques-Severi-Zariski lemma  (see \cite[III.7.8]{hartsh} or \cite[\href{https://stacks.math.columbia.edu/tag/0FD9}{Tag 0FD9}]{stacks-project}).

Theorem~\ref{chev.thm}  was proved by Chevalley  and Barsotti around 1953, with proofs appearing in
\cite{MR76427, MR82183, MR126447, MR159825}, using classical terminology and methods.
Some of the history is discussed in \cite[VII.5]{MR1847105}, and   modern versions are  in \cite{MR1906417, MR3088271, milne2013proofbarsottichevalleytheoremalgebraic}.

The existence of Jacobians
over $\c$  goes back to Abel and Jacobi. See \cite{MR2077579} for a historical essay, and \cite[Chap.4]{MR0257326} or
  \cite[Sec.2.7]{gri-har} for proofs.  Over arbitrary fields the construction is due to Weil \cite{MR29522}. 
The general case is stated,
but not  fully treated, in  \cite{shaf, hartsh}.

\end{say}

\section{Algebraic geometry background}

We need only one result form algebraic geometry that is not treated in
basic textbooks:
 the existence of the Jacobian, which is the group of degree 0 divisors modulo linear equivalence.

 \begin{thm}\label{jac.exists}
   Let $X$ be a  smooth, projective curve. Then  
  its Jacobian is  a proper, algebraic group, denoted by $\jac(C)$. \qed
 \end{thm}

 The proof in \cite{MR29522} is very clear; so is the  version in  \cite{MR861976} which uses contemporary language.
It was superseded by the construction of the Picard variety in \cite{fga}; this viewpoint is presented in \cite{mumf66, blr}.

 \begin{notation}\label{notation.1}
   Let $k$ be a field with algebraic closure $\bar k$.
   
   A {\it $k$-variety} is an irreducible and reduced scheme $X$ of finite type over $k$.  Its {\it smooth locus} is denoted by $X^{\rm sm}$.

   The adjective {\it geometrically} refers to properties over $\bar k$.

   A {\it rational map} between varieties is denoted by a dotted arrow
   $\phi:X\map Y$.  A {\it morphism} is everywhere defined, and  denoted by a solid arrow
   $\phi: X\to Y$.

   The {\it fiber}  $\phi^{-1}(y)$ over $y\in Y$ is usually denoted by $X_y$.

   The notation for {\it linear systems} is discussed in Section~\ref{weil.div.sect}.

   Let $S$ be a $k$-scheme and  ${\mathcal P}$ a property that makes sense for $K$-points of $S$. 
 We say that  ${\mathcal P}$ holds for {\it general points} of $S$, iff  there is a dense,  open  subset  $S^\circ\subset S$ such that
  ${\mathcal P}$ holds  for  all $K$-points of $S^\circ$ for every field extension $K\supset k$.

   Let $G$ be a connected, smooth, algebraic group with unit $e\in G$.
   A {\it rational action} on a normal variety $X$ is a rational map
   $$   \mu: G\times X\map X$$ such that  $\mu(e,x)=x$ and 
   $\mu(g_1g_2, x)=\mu\bigl(g_1, \mu(g_2, x)\bigr)$
   holds for $(g_1, g_2, x)$ in a dense, open subset of $G\times G\times X$.  One can identify $\mu$ with the collection of birational maps
   $$
   \mu_g:X\cong \{g\}\times X \stackrel{\mu}{\map} X;
   $$
these are defined for $g$  in a dense, open subset of $G(\bar k)$.

   In classical terminology a morphism is called a `regular map.' I try to avoid this, since both \cite[IV.6.8.1]{ega} and \cite[\href{https://stacks.math.columbia.edu/tag/07R6}{Tag 07R6}]{stacks-project} use `regular' with a very different meaning.
However, if $\mu$ is a morphism, it is usually called an {\it action} 
or {\it regular action} to emphasize the difference from a rational action.

Let $H, G$ be algebraic groups and $\phi:H\map G$ a rational map.
It is a {\it rational homomorphism} of algebraic groups if there are  dense, open subsets
  $U\subset H$ and   $W\subset U\times U$, 
such that $\phi(u_1u_2)=\phi(u_1)\phi(u_2)$ for $(u_1, u_2)\in W$.
If $\phi$ is morphism, then we get the usual notion of 
homomorphism or {\it  regular homomorphism.} 
(We see in Lemma~\ref{rtl.h.h.lem} that a rational group homomorphism is  actually regular, but in Lemma~\ref{lin.rep.cor} we will obtain a map
$H\map \PGL_N$ which is only a rational homomorphism to start with.)
   \end{notation}

 \section{Proof of Chevalley's Theorem}

 \begin{notation} We work over  a perfect field  $k$ with algebraic closure $\bar k$. Let $G$ be a connected, smooth, algebraic group, with unit $e\in G$. The group operation is written multiplicatively.

   A proper algebraic group is frequently denoted by $A$, and the group operation is written additively. In the proof these will always appear as
   subgroups of products of Jacobians of smooth, projective curves.
 \end{notation}

 \begin{say}[Proof of Theorem~\ref{chev.thm}] \label{chev.thm.pf}
The proof  is done in 5 steps.
 
\medskip
{\it Step \ref{chev.thm.pf}.1: Rational action on a projective  variety.}
     The arguments are the simplest if we have a regular, faithful $G$-action on a normal, projective variety $X$. However, if we are willing to work with Weil divisors, then it is enough to use a {\em rational,}
   faithful $G$-action on a normal, projective variety $X$. These are easy to get.

 Let $U\subset G$ be any affine, open subset, and  $X$ the normalization of its closure in some projective space.  Then $X$ is normal and projective, with a  rational,  faithful $G$-action.

\medskip
{\it Step \ref{chev.thm.pf}.2: Algebraic map of $G$ to an Abelian variety.}
   Let $G$ be a connected, smooth, algebraic group acting rationally on a normal, projective variety $X$. We write the action as
   $\mu: G\times X\map X$.

   Let $D$ be a Weil divisor on $X$. The closure of the image of $G\times D$ by $(1_G, \mu)$ is a  Weil
   divisor ${\mathbf D}$ on $G\times X$. Let ${\mathbf D}_g$ denote its restriction to $\{g\}\times X\cong X$; see Definition~\ref{wfam.defn}.

   For a smooth, projective curve $C\subset X^{\rm sm}$, let
   $\alpha_C: G\to \jac(C)$ be the morphism that sends
  $g\in G$ to  ${\mathbf D}_g|_{C}-{\mathbf D}_e|_{C}$.

  The Severi-Weil  Lemma~\ref{weak.w-m.lem} says that there is a finite set of curves
    $\{C_i\subset X^{\rm sm}: i\in I\}$ with  product morphism
   $$
     \alpha_{X, D}: G\to  \tprod_{i\in I}\jac(C_i),
     \eqno{(\ref{chev.thm.pf}.2.1)}
     $$
 such that $\alpha_{X, D}(g_1)=\alpha_{X, D}(g_2)$ iff   $ {\mathbf D}_{g_1}$ is linearly equivalent to ${\mathbf D}_{g_2}$.
 Note that  $\alpha_{X, D}(e)=0$.

 At this point $\alpha_{X, D}$ is a morphism of varieties, and we do not yet claim  any naturality property for it.

\medskip
{\it Step \ref{chev.thm.pf}.3:  Homomorphism of $G$ to an Abelian variety.}
   The morphism  $\alpha_{X, D}$  obtained in (\ref{chev.thm.pf}.2.1) is 
   a group homomorphism, since any morphism from an  algebraic group to a proper, algebraic group is a group homomorphism by 
 \cite{MR29522}. We recall the proof in Section~\ref{sec.p1.pf}.

\smallskip

Let $H=H_{X, D}\subset G$ denote the
identity component of the kernel of $\alpha_{X, D}$.
We aim to show that $H$ is affine for suitable choice of $D$.

\medskip
    {\it Step \ref{chev.thm.pf}.4: Linear representations of $H$.}

    We start with a special case when
   the $G$-action on $X$ is regular. Then translation by any $h\in G$
   gives  isomorphisms of the complete linear systems
   $$
   \mu_{h,g,D}: |{\mathbf D}_g|\stackrel{\cong}{\longrightarrow} |{\mathbf D}_{hg}|.
   $$
   If  $h\in H$ then  $|{\mathbf D}_g|=|{\mathbf D}_{hg}|$, so for any fixed $g$ we get a
   linear representation
   $$
   \tau_{g,D}: H\to  \aut\bigl(|{\mathbf D}_g|\bigr)\cong \PGL_N.
$$
    The elements of the kernel  map $D$ to itself.
A sufficiently general $D$ should have no automorphisms. If this holds then  $\tau_{g, D}$ is faithful, hence $H$ is affine.

    \smallskip

    If  the $G$-action on $X$ is not regular, then the best 
    we can hope for is that the $\mu_{h,g,D}$  are  isomorphisms 
    for $(h, g)$ in some open, dense subset of $H\times G$.
    However, even this seems problematic, because of the following subtle  difficulty.
    \smallskip
    
 {\it Warning \ref{chev.thm.pf}.4.1.} 
Let $\{A_i: i\in I\}$ be pairwise linearly equivalent divisors spanning the  linear system  $|A_i: i\in I|$. Even if the $g_*(A_i)$ are  pairwise linearly equivalent divisors, the linear system spanned by the $g_*(A_i)$ is frequently different from the birational transform of the  linear system  $|A_i: i\in I|$. That is, 
   $$
g_*|A_i: i\in I|\neq |g_*(A_i): i\in I|,
\eqno{(\ref{chev.thm.pf}.4.2)}
$$
see Example~\ref{bir.tr.weil.2.ex}.
However, Corollary~\ref{bir.tr.weil.2} shows that equality holds in
(\ref{chev.thm.pf}.4.2) if the $A_i$ are parametrized by a geometrically  irreducible variety.

    A solution is to replace  the complete linear system $|{\mathbf D}_g|$ with a  carefully  chosen linear subsystem.

\medskip
{\it Step \ref{chev.thm.pf}.4.3: Linear representations of $H$, general case.}
   With $G\supset U\subset X$ as in (\ref{chev.thm.pf}.1),
   there is a dense, open subset $G^*\subset G$ such that
   ${\mathbf D}_g$ is the closure of its intersection with $U$
   for $g\in G^*$. Thus, if $g, hg\in G^*$, then
   $(\mu_h)_* {\mathbf D}_g={\mathbf D}_{hg}$.
   However, as we noted in (\ref{chev.thm.pf}.4.1), this does not imply that $(\mu_h)_* |{\mathbf D}_g|$ and $|{\mathbf D}_{hg}|$ are the same.

   Fix now a $g\in G^*$. Let $H^*\subset H$ be a dense, open subset such that $H^*\cdot g\subset G^*$.
Let $H^{**}\subset H^*$ be any dense, open subset.
   Then the linear subsystem
   $$
   |{\mathbf D}_g|^H:=\bigl| {\mathbf D}_{hg}: h\in  H^{**}\bigr|
   \subset |{\mathbf D}_g|
   \eqno{(\ref{chev.thm.pf}.4.4)}
   $$
   is independent of $H^{**}$.
   In particular, $\mu_{h}$ maps a  spanning set parametrized by
   $H^*\cap h^{-1}H^*$ to the  spanning set parametrized by
   $hH^*\cap H^*$.
Since $H$ is geometrically  irreducible, 
   Corollary~\ref{bir.tr.weil.2}  implies that 
    $$
   (\mu_h)_* |{\mathbf D}_g|^H=|{\mathbf D}_g|^H \qtq{for} h\in H^*.
   \eqno{(\ref{chev.thm.pf}.4.5)}
   $$

   Once   (\ref{chev.thm.pf}.4.5) holds, Corollary~\ref{lin.rep.cor} says that we get a 
   linear representation
   $$
   \tau_{g,D}: H\to  \aut\bigl(|{\mathbf D}_g|^H\bigr)\cong \PGL_N.
     \eqno{(\ref{chev.thm.pf}.4.6)}
$$
 Its   kernel is a subgroup  of those elements of $H$ that map $D$ to itself.

\medskip
{\it Step \ref{chev.thm.pf}.5: Choosing $D$.}
   For `general'  $D$ we expect that  the kernel $K$ of
   $\tau_{g,D}$ is trivial, but here is a concrete example.

   Let $U\subset G$ be an affine, open subset as in (\ref{chev.thm.pf}.1).  Choose $D_U\subset U$ to have an isolated singular point, and let $D\subset G$ be  its closure.

   If a subgroup $K\subset H$ fixes $D$, then it fixes $\sing D$.
    So, if $\sing D$ has   an isolated point, then $K$ is finite.
    Then $\tau_{g,D}$ is  a linear representation of $H$ with a finite kernel, hence $H$ is affine.
     This completes the proof of Theorem~\ref{chev.thm}. \qed
\smallskip
    
As a slight improvement, 
    if $\chr k\neq 2$ or $\chr k= 2$ and $\dim G$ is even, then
    we can chose $D_U$ such that its scheme theoretic singular set is a closed point. In this case we also get that $K$ is reduced.
    If $\chr=2$,  we can replace $G$ by $G\times \gm$.
 \end{say}

 \begin{rem}\label{chev.thm.rem.1} The stronger form of Theorem~\ref{chev.thm} says that
   $G$ has a unique, connected, affine, normal  subgroup $H$ such that $G/H$ is proper.

   Our proof shows that there is a minimal such $H^{\rm min}\subset G$.
   This implies the stronger form, if we use
   that the quotient of an affine, algebraic group by a normal, algebraic  subgroup
   is again an affine, algebraic group; see 
   \cite[Sec.11.5]{humhreys} or  \cite[Thm.6.8]{borel}. 
 \end{rem}

 \section{Maps to  algebraic groups}\label{sec.p1.pf}

 These results are mostly  due to Weil \cite{MR29522}; see also \cite[Sec.4.4]{blr}. 
The proofs are reproduced
 for completeness.

  \begin{lem}\label{abv.lem} Let $(G, \cdot, e)$ be a connected, algebraic group, $(A, +, 0)$ a proper, algebraic group and $\phi:G\to A$ a morphism such that $\phi(e)=0$. 
    Then $\phi$ is a group homomorphism.
     \end{lem}

  Proof. 
  Let  $\tau:G\times G\to A$ be the map
  $$
  \tau:(g_1,g_2)\mapsto \phi(g_1g_2)-\phi(g_2)-\phi(g_1).
  $$
  We  need to show that $\tau$ is constant on a dense, open subset of $G\times G$.  Choose an affine $e\in U\subset G$ and a normal, projective
  $U\into \bar U$. Given $g_1,g_2\in U$ choose
  smooth, projective, irreducible curves 
  $\{e,g_i\}\subset C_i\subset \bar U$.
Then $\tau$ restricts to a rational map
$\tau_C:C_1\times C_2\map A$, which is a morphism by (\ref{weil.morph.lem}).

Also, $\tau_C$ maps  both 
$\{e\}\times C_2$ and $C_1\times \{e\}$ to $0\in A$.
We apply the Rigidity lemma  of \cite[Sec.4]{mumf-abvar} and \cite[Sec.III.4.3]{shaf}. We get that 
$\tau_C$ maps both 
$\{c_1\}\times C_2$ and $C_1\times \{c_2\}$ to points for all $c_i$.
Thus $\tau_C$ is constant. \qed

\begin{lem}\label{weil.morph.lem} Let $X$ be a smooth variety and $A$ a proper algebraic group. Then every rational map $\phi:X\map A$ is a morphism.
\end{lem}

Proof. Since $A$ is proper, there is a closed subset $Z\subset X$ of codimension $\geq 2$ such that $\phi|_{X\setminus Z}$ is a morphism.
Let $\Delta_X\subset X\times X$ denote the diagonal, and 
consider
$$
\Phi:X\times X\map A \qtq{given by} \Phi(x_1, x_2):=\phi(x_1)-\phi(x_2).
$$
Note that $\Phi$ is a morphism on $(X\setminus Z)\times (X\setminus Z)$
and maps $\Delta_{X\setminus Z}$ to $0\in A$.

Let $0\notin D\subset A$ be a Cartier divisor such that
$A\setminus D$ is affine. Let $\Phi^*(D)\subset X\times X$ denote the closure of its preimage.  Set $W:=(X\times X)\setminus \Phi^*(D)$.
Then $\Phi|_W$ is a rational map to the affine variety $A\setminus D$,
and it is a morphism  outside the codimension $\geq 2$ set
$(X\times Z)\cup (Z\times X)$. Thus $\Phi|_W$ is a morphism.

Since $X\times X$ is smooth, $\Phi^*(D)$ is a Cartier divisor.
Thus $\Phi^*(D)\cap \Delta_X$ has pure codimension 1 in $\Delta_X$.
On the other hand, $\Phi^*(D)\cap \Delta_X$
is contained in  $\Delta_Z$ which has codimension  $\geq 2$ in $\Delta_X$. Thus $\Phi^*(D)$ is disjoint from $\Delta_X$.
So $\Phi$ is a morphism in a neighborhood of $\Delta_X$. 

Finally note that 
$\phi(x)=\Phi(x, y)+\phi(y)$, and for every $x$ there is a dense, open subset
$U_x\subset X$ such that $\Phi(x, y)$ and $\phi(y)$ are both defined for $y\in U_x$.
So $\phi$ is a morphism. \qed

\medskip
{\it Remark  \ref{weil.morph.lem}.1.} The proof gives the following stronger version.

  Let $X$ be a normal variety such that $X\times X$ is $\q$-factorial, and $A$ an algebraic group. Let $\phi:X\map A$ be a  rational map
    that is  a morphism outside a  codimension $\geq 2$ subset. Then  $\phi$ is a morphism everywhere.

\begin{lem}[Rational homomorphisms]\label{rtl.h.h.lem}
  Any rational homomorphism of algebraic groups is a regular  homomorphism.
\end{lem}

Proof. 
Let  $\phi:H\map G$ a rational homomorphism as in Notation~\ref{notation.1}.
  Consider the rational map
  $$
  \Psi: (u_1, u_2, u_3)\mapsto \phi(u_1u_3)\phi(u_3^{-1}u_2).
  $$
  It is independent of $u_3$ on a dense, open subset, hence everywhere.
  Thus $\Psi$  depends only on the product $u_1u_2$, and gives the extension of $\phi$ to $H$. \qed

  \begin{cor}[Linear representations from rational actions] \label{lin.rep.cor}
    Let $\mu:G\times X\map X$  be a connected algebraic group acting rationally on a  normal, projective variety.
    Let $|D|$ be a linear system on $X$ such that
    $(\mu_g)_*|D|=|D|$ for general $g\in G(\bar k)$. Then
    $g\mapsto (\mu_g)_*$ defines a linear representation
    $$
    \tau_{|D|}:  G\to \aut(|D|)\cong \PGL_N,\qtq{where} N:=\dim |D|+1.
    $$
  \end{cor}

  Proof. Fix two  point sets $\{p_i\in \p^N:i=1,\dots, N+2\}$ and
    $\{q_i\in \p^N:i=1,\dots, N+2\}$  in
  linearly general position. Then there is a unique $\phi\in \aut(\p^N)$ such that $\phi(p_i)=q_i$ for every $i$.

  Thus giving a rational map $\Phi: G\map  \aut(\p^N)$ is the same as giving $N+2$ rational maps $\phi_i:G\map \p^N$ such that
  the $\{\phi_i(g):i=1,\dots, N+2\}$ are in
  linearly general position  over a dense, open subset $G^*\subset G$.

  By assumption we have a rational map
  $$
  \bar\mu: G\times |D|\map |D|.
  $$
  We may assume that $\bar\mu$ is defined on a dense, open subset of
  $\{e\}\times |D|^\circ \subset \{e\}\times |D|$. Choose $N+2$ divisors $D_i\in |D|^\circ$ in linearly general position. Then the
  $S_i:=\bar\mu\bigl(G\times [D_i]\bigr)$ are $N+2$ rational sections
  of $\pi_1:G\times |D|\to G$.

  There is thus a dense, open subset $G^*\subset G$ such that
  the $S_i$ are sections in  linearly general position over $G^*$.
  Thus the $S_i$ define a  rational representation
  $\tau_{|D|}:  G\map \aut(|D|)$, which is   regular 
  by Lemma~\ref{rtl.h.h.lem}. \qed

 \section{Linear systems of  Weil divisors}\label{weil.div.sect}

  \begin{defn}[Weil divisors]\label{wfam.defn}  
   Let $X$ be a normal variety. A {\it Weil divisor} on $X$ is a finite $\z$-linear combination of irreducible,  codimension 1 subvarieties   $D:=\sum m_i D_i$.

   Let $X^{\rm sm}\subset X$ be the smooth locus. Then $D|_{X^{\rm sm}}$ is a Cartier divisor, whose closure is $D$. One can thus identify Weil divisors on $X$ with Cartier divisors on $X^{\rm sm}$; though this can be misleading for families of divisors.

   A {\it family of Weil divisors}  parametrized by a smooth variety $S$ is a Weil divisor ${\mathbf D}$ on $X\times S$ whose support does not contain any of the fibers $X_s:=X\times \{s\}$.  (This definition also works if $S$ is normal, but needs more foundational work, see \cite[Sec.4.3]{k-modbook}.)
 The fiber ${\mathbf D}_s$ over $s\in S$ is defined as the the closure of the restriction of ${\mathbf D}$ to 
 $X_s^{\rm sm}$.

 Let $U\subset X\times S$ be any dense, open subset.
   Then ${\mathbf D}_s\setminus U$ has codimension $\geq 1$  for general $s$. Thus

   \smallskip
    {\it Claim \ref{wfam.defn}.1.}
${\mathbf D}_s$ is the closure of ${\mathbf D}_s\cap U$  for all $s$ in  some dense, open $S^\circ\subset S$.
\end{defn}

 \begin{defn}[Pull-back and  birational transform]\label{w.birtr.defn}  
   Let $h:Y\to X$ be a morphism of normal varieties.
   Assume that there is a closed  $Z\subset X$ such that
   $D$ is  Cartier on $X\setminus Z$ and 
   $h^{-1}(Z)$ has codimension $\geq 2$ in  $Y$. Then one can define the {\it pull-back}  $h^*(D)$ as the closure of the Cartier divisor
   $h^*(D|_{X\setminus Z})$.    If these condition do not hold, then one usually cannot define
   the pull-back of $D$ in a sensible way.

   The pull-back preserves linear equivalence.

 Let $X, Y$ be normal varieties. Let $g:X\map Y$  be a birational map that  is defined outside a closed subset $Z\subset X$ of codimension $\geq 2$. (For example, $Y$ is proper.)
    If $D$ is a Weil divisor on $X$, we define its 
    {\it birational} or {\it proper} transform on $Y$ as
    the closure of $g_*(D|_{X\setminus Z})$. We denote it by
    $g_*(D)$.

      \smallskip
 {\it Warning.}  Birational transform  need not   preserve linear equivalence. So we need  to define the birational transform of a linear system separately in (\ref{bir.tr.weil}). 
 \end{defn}

 \begin{defn}[Linear systems]\label{lin.sys.defn}
   Let $X$ be a normal, proper and geometrically integral $k$-variety.
We let $|D|$ denote a not-necessarily complete linear system.
   Its members are usually denoted by   $D_\lambda$ as divisors, and
   $[D_\lambda]$ as points in the projective space $|D|$.

   \medskip

   {\it Warning.} In the literature $|D|$ is used to denote either the {\it complete} linear system of a divisor $D$, or some not-necessarily complete linear system; in many cases without explicit mention.
   I intend to say complete   linear system when that variant is used.
    \medskip
   
   If $\{D_\lambda:\lambda\in \Lambda\}$ is a set of linearly equivalent Weil divisors, we let $|D_\lambda:\lambda\in \Lambda|$ denote the
   (not-necessarily complete) linear system spanned by them.

   The {\it incidence correspondence}  of $|D|$ is
 $$
 {\mathbf D}:=\{(x, D): x\in D\}\subset X\times |D|.
 $$
  If $|D|$ has no fixed components
then ${\mathbf D}$ is irreducible and its
  first projection  $\pi_2: {\mathbf D} \to X$
  is generically a hyperplane bundle.
  Conversely, any such ${\mathbf D}$ is the incidence correspondence of a unique linear system $|D|$ without fixed components.

 A 1-dimensional linear system is traditionally called a {\it pencil.}
 \end{defn}

 \begin{say}[Birational transforms of linear systems]\label{bir.tr.weil}
   Let  $g:X\map Y$ be a birational map between
normal, proper and geometrically integral $k$-varieties,
and    $|D|$ a  linear system on $X$ without fixed components.
    
    Using the identity of $|D|$, we get a birational map
    $\bar g:X\times |D|\map Y\times |D|$.  Then  $\bar g_*({\mathbf D})$
    is a divisor on  $Y\times |D|$, which corresponds to a  linear system  $g_*|D|$ on $Y$. It is   called the {\it birational} or {\it proper} transform of $|D|$.
 \end{say}

The next example illustrates some of the subtleties that one needs to keep in mind.

 \begin{exmp}\label{bir.tr.weil.2.ex} Let
   $\phi:(x{:}y{:}z)\mapsto (x^{-1}{:}y^{-1}{:}z^{-1})$ be the standard quadratic transformation of $\p^2$. The  linear system of lines
   $|L|:=|ax+by+cz|$ is mapped to the linear system of conics
   $|Q|:=|ayz+bxz+cxy|$. Let $W$ be the set of lines that pass through a coordinate vertex.  It is a union of 3 linear subsystems
   $|ax+by|\cup|by+cz|\cup |ax+cz|$.
   Each  of these is mapped to  itself, for example
   $|ax+by|$ is mapped to $|ay+bx|$.

   Thus $W\subset |L|$ is a connected subset that spans $|L|$, yet the $\phi$-images of its members do not span $|Q|=\phi_*|L|$.

   If we change $\phi$ to $\phi'$ so that its 3 base points are conjugate over $k$, then the corresponding $W'$ is irreducible over $k$. 
   \end{exmp}

 In working with linear systems, we need to understand
 when the birational transform of a member of a linear system $|D|$ fails to be a
 member of the birational transform of  $|D|$.
 The next lemma says that this happens only along linear subspaces.
 As an immediate corollary we get that
general members of 
 irreducible spanning sets transform as expected.

    \begin{lem}\label{bir.tr.weil.1}
 Let  $g:X\map Y$ be a birational map between normal $k$-varieties as in (\ref{lin.sys.defn}), and 
   $|D|$ a    linear system on $X$ without fixed components.
          Then there is a finite set of linear subspaces $\{L_i\subsetneq |D|_{\bar k}: i\in I\}$ defined over $\bar k$, such that  $g_*(D_\lambda)\in g_*|D|$ 
  for $[D_\lambda]\notin  \cup_i L_i$. 
  \end{lem}

 Proof. The map $g$ restricts to an isomorphism on suitable,
 open, dense subsets  $X\supset U_X\cong U_Y\subset Y$.
  Let $E_i$ be the codimension 1, irreducible components of
  $Y\setminus U_Y$, defined over $\bar k$.
  
   Then
  $|D|$ and $g_*|D|$ restrict to the same set of divisors on
   $U_X$ and $U_Y$.  Thus if $[D^X_\lambda]\in |D|$, then
   the corresponding member of $[D^Y_\lambda]\in  g_*|D|$ is of the form
   $$
   D^Y_\lambda=g_*\bigl(D^X_\lambda\bigr)+\tsum m_{i,\lambda} E_i.
   $$
   Let
   $L^Y_i\subset g_*|D|_{\bar k}$ be the linear subspace of those divisors that contain $E_i$.  Since $g_*|D|_{\bar k}$ has no fixed components, $L^Y_i\neq g_*|D|_{\bar k}$.
   We can then take $L_i$ to be the corresponding linear subspace of $|D|_{\bar k}$. \qed

 \begin{cor}\label{bir.tr.weil.2}
Let $|D|$ be a  linear system on $X$ without fixed components, and $g:X\map Y$ a birational map. Let $W\subset |D|$ be a geometrically  irreducible subvariety that spans $|D|$.
Then there is a dense, open $W^\circ\subset W$ such that
$$
g_*|D|=\bigl|g_*(D_\lambda): [D_\lambda]\in W^\circ\bigr|.
\qed
$$
 \end{cor}

 \section{Bertini theorems}\label{bert.sec}

 Let $X$ be a variety and $|H|$ a linear system on  $X$.
 Bertini-type theorems say that if $X$ satisfies certain property ${\mathcal P}$, and $|H|$ is `large enough,' then `general members'  $[H_{\lambda}]\in|H|$ also satisfy ${\mathcal P}$. Working over a field $k$, we usually want the $H_{\lambda}$ to be defined over $k$.
This works whenever $k$ is  infinite.
For Bertini-type theorems over finite fields, see \cite{cha-poo}.

We will use the case ${\mathcal P}=(\mbox{smooth})$, which   is proved in many places, for example \cite[Sec.II.6.2]{shaf}. 
\cite[II.8.18]{hartsh} states it for $X$ projective, but the proof is the same in general; we essentially give it in (\ref{bert.irr.thm.cor}). (Note that, in the quasi-projective case, the set  $\{H_{\lambda}: X\cap H_{\lambda} \mbox{ is smooth}\}$ need not be open in  $|H|$.)

\begin{thm}[Bertini's smoothness theorem]\label{bert.sm.thm}
   Let $X\subset \p^N$ be a smooth, quasi-projective variety and
    $H\subset \p^N$   a general hyperplane. Then
   $X\cap H$ is smooth. \qed
\end{thm}

Similar results for other local properties ${\mathcal P}$  are treated in  \cite[Sec.IV.12.1]{ega}; see also \cite[Sec.10.2]{k-modbook}.

 We start the proof of Bertini's irreducibility theorem with some lemmas.

\begin{lem}\label{sm.p.g.i.lem}  Let $Y$ be an irreducible $k$-variety that has a smooth $k$-point $p$. Then $Y$ is geometrically irreducible.
\end{lem}

Proof. Let $k^{\rm s}$ be the separable closure and $Y^{\rm s}_i$ the irreducible components of $Y_{k^{\rm s}}$. One of them, say $Y^{\rm s}_1$ contains $p$. Then all of its conjugates also contain $p$. Since $p$ is smooth, $Y^{\rm s}_1$ is
$\gal(k^{\rm s}/k)$-invariant, hence defined over $k$. So
$Y^{\rm s}=Y^{\rm s}_1$. \qed

\begin{lem}\label{bert.irr.lem}
  Let $X$ be a geometrically irreducible variety and
  $|M|$ a pencil  of divisors without fixed components. Assume that members of $|M|$ are smooth at a $k$-point $x\in X$, and meet transversally. (That is, the scheme-theoretic  base locus $\bs|M|$  is smooth  and has codimension 2  at $x$.)
  Then general members of $|M|$ are geometrically irreducible.
  \end{lem}

Proof. We may assume that the base field is algebraically closed.
By shrinking $X$ we may assume that $\bs|M|$ has pure codimension 2 in $X$.

Let $\bar X\subset X\times \p^1$ be the closed graph of $X$.
The first projection $\pi_1:\bar X\to X$ is an isomorphism over
$X\setminus \bs|M|$, and a $\p^1$-bundle over $\bs|M|$.
Thus $\bar X$ is (geometrically) irreducible.
Therefore the generic fiber $\bar X_\eta$ of the second projection $\pi_2:\bar X\to \p^1$ is
irreducible.  Finally $\pi_1^{-1}(x)$ gives a smooth $k(\eta)$-point on $\bar X_\eta$. Thus 
 $\bar X_\eta$ is geometrically irreducible by Lemma~\ref{sm.p.g.i.lem}.

Then general members of $|M|$ are geometrically irreducible, see for example  \cite[Sec.II.6.1]{shaf}.
\qed

\begin{say}[Proof of  Theorem~\ref{bert.irr.thm}]\label{bert.irr.thm.pf}
   We may assume that the base field is algebraically closed.
By shrinking $X$ we may also assume that each fiber of $\pi$ has pure codimension $n$ in $X$.

Assume first that $\pi$ is generically smooth, and choose  $p\in \p^n$ such that $X_p$ is generically smooth. 
Let $|M|\subset |\o_{\p^n}(1)|(-p)$ be a general pencil and
$|M|_X$ its pull-back to $X$. Then  members of $|M|_X$ are smooth at $x$ and meet transversally. Thus general members of $|M|_X$ are
irreducible by (\ref{bert.irr.lem}).

If $\pi$ is not generically smooth, then consider the $p^r$-th power  Frobenius
$F^r:\p^n_{(r)}\to \p^n$.  (In fact $\p^n_{(r)}\cong \p^n$.)  Let $X_{(r)}$ denote the reduced underlying subscheme of
$\p^n_{(r)}\times_{\p^n} X$. Then $X_{(r)}\to X$ is a homeomorphism,
so $\pi^{-1}(H)$ is geometrically irreducible iff its pull-back to $X_{(r)}$ is. Finally 
$X_{(r)}\to\p^n_{(r)}$ is  generically smooth for $r\gg 1$, so the previous argument applies. \qed
 \end{say}

Putting these together and using induction we get that
 if $X\subset \p^N$ is a normal, projective variety of dimension $n$,  and
 $H_1, \dots, H_{n-1}\subset \p^N$   are general hyperplanes,
 then $X\cap H_1\cap \dots\cap H_{n-1}$ is a smooth, irreducible curve.

 We need  a stronger version of this, where the $H_i$ are not `fully general.' To state it,
 let $X$ be a  projective variety of dimension $n$, and
$|H|$ a  linear system.
Let $\grass(n{-}1, |H|)$ denote the Grassmannian of $n{-}1$ dimensional subspaces of $|H|$. Let ${\mathbf C}\subset X\times \grass(n{-}1, |H|)$ denote the  universal subvariety of curve intersections. That is, the fiber $C_L$ of  ${\mathbf C}$ over $L\in \grass(n{-}1, |H|)$ is the intersection of all members of $L$.

\begin{cor}[Bertini's irreducibility theorem II]\label{bert.irr.thm.cor}
  Let $Y$ be a  projective variety  such that $\sing X\subset X$ has codimension $\geq 2$. Let $A$ be a  very ample  divisor and $|H|:=|mA|$, the complete linear system. 
  Then there are closed subsets  $Z_i\subset \grass(n{-}1, |H|)$ of codimension $\geq i$ such that
  \begin{enumerate}
  \item $C_L$ is smooth and irreducible if $m\geq 1$ and $L\notin Z_1$,
     \item $C_L$ is irreducible if  $m\geq  2$ and $L\notin Z_2$, except for $Y=\p^2$ and $|H|=|\o_{\p^2}(2)|$.
\end{enumerate}
  \end{cor}

Proof. It is enough to show that if $S\subset Y$ is the intersection of
$n{-}2$ general members of $|H|$ then the claims hold for $S$.
Note that $S$ is irreducible by Theorem~\ref{bert.irr.thm}, and has only
isolated singularities by (\ref{bert.sm.thm}).

Let  $|H_S|$ be the restriction of  $|H|$ to $S$.
Using Theorem~\ref{bert.irr.thm} on $S$, we see that  if 
$C\in |H_S|$ is reducible, then $C$ either passes through a singular point of $S$, or it is singular at a smooth point of $S$.

If $s\in \sing S$ then $|H_S|(-s)$ is a hyperplane in  $|H_S|$, and a general member of $|H_S|(-s)$ is irreducible, unless
$|H_S|(-s)$ maps $S$ to a  curve. So $S$ is a cone with vertex $s$, and 
 $H$ has degree 1 on the rulings. This cannot happen for $m\geq 2$.

If $s\in S\setminus\sing S$ then $|H_S|(-2s)$ has codimension 3 in  $|H_S|$.
A general member of $|H_S|(-2s)$ is irreducible, unless
 $S$ is a cone with vertex $s$, and 
 $H$ has degree $\leq 2$ on the rulings.
This cannot happen for $m\geq 3$.

For  $m=2$  the only problem is when
$S$ is a cone with vertex $s$
 for every $s$. This happens only for $S\cong\p^2$. \qed

 \section{Severi-Weil lemmas}\label{alb-pic.sect}

  Severi noted that linear equivalence of divisors on surfaces can be checked on finitely many   curves  \cite{sev-29, sev-35, sev-94}. This was extended by Weil to
  normal varieties \cite{MR47368}.
  
  Note that while the restriction of a Weil divisor $D$ to a curve $C$ is  usually not defined, $D$ is Cartier on the smooth locus $X^{\rm sm}$,
  thus $D|_C$ is a well-defined linear equivalence class whenever
  $C\subset X^{\rm sm}$.

Given a  morphism $\pi:Y\to S$, we would like to know when a divisor is 
linearly equivalent  to the pull-back of a divisor on $S$.
As first step, we look for   divisors that do not dominate $S$.
Equivalently, they are disjoint from the generic fiber.
I call these {\it vertical.}

 \begin{lem} \label{ag.lem.1} Let $Y, S$ be  normal varieties,  and $\pi:Y\to S$ a proper morphism whose general fibers are geometrically integral curves.  Let $D$ be a Weil divisor on $Y$. Then $D$ is linearly equivalent to a vertical divisor iff
   $D|_{Y_s}\sim 0$  for general $s\in S$.
 \end{lem}

 Proof. If $D$ is vertical, then it is disjoint from a general fiber $Y_s$, so $D|_{Y_s}=0$. Conversely, if $D|_{Y_s}\sim 0$  for general $s\in S$,
 then $D|_{Y_g}\sim 0$  for the generic fiber. 
  This gives a rational section of $\o_Y(D)$ whose divisor consists of $D$ plus a  vertical divisor. \qed

 \begin{lem} \label{ag.lem.2} Let $Y, S$ be  normal varieties,  and $\pi:Y\to S$ proper morphism of relative dimension 1.
   Assume that there is a closed subset  $Z_2\subset S$ of  codimension $\geq 2$,  such that  the fibers of $\pi$ are
   geometrically integral over $S\setminus Z_2$.

Let $D$ be a Weil divisor on $Y$. Then $D$ is linearly equivalent to the
pull-back of a divisor from $S$ iff
   $D|_{Y_s}\sim 0$  for general $s\in S$.
 \end{lem}

 Proof. The assumption implies that every  vertical divisor is the 
 pull-back of a divisor from $S$. \qed

 \medskip

 \begin{lem}\label{weak.cl0.prop}
   Let $X$ be a geometrically normal, projective variety over an infinite field $k$, and  $|H|$ a very ample linear system on $X$.
   Let  $D$ be a  Weil divisor on $X$. The following are equivalent.
   \begin{enumerate}
   \item $D\sim 0$.
     \item $D|_C\sim 0$ for every 
smooth, complete intersection curve
$C:=H_1\cap\cdots \cap H_{n-1}$ over $k$, where $H_i\in |2H|$.
\end{enumerate}
   \end{lem}

   Proof. It is clear that (\ref{weak.cl0.prop}.1) implies
   (\ref{weak.cl0.prop}.2), and the equivalence is clear 
 for $X=\p^2$.

   For the converse $|2H|$ gives an embedding  $X\into \p^N$. Choose a
   general projection $\pi:X\map \p^{n-1}$ whose base point set $P=\{x_i\}$ is smooth.
   We may assume that $P$ is disjoint from $\supp D$.
   Blowing up $P$ we get $X':=B_PX$ and $D'\cong D$. 

The projection  morphism
   $\pi':X'\to \p^{n-1}$ has relative dimension 1,
   and satisfies the assumptions of Lemma~\ref{ag.lem.2} by Corollary~\ref{bert.irr.thm.cor}.
     Thus, if (\ref{weak.cl0.prop}.2) holds, then $D'$  is
   linearly equivalent to the
   pull-back of a divisor $D_P$ from $\p^{n-1}$.

   Over each $x_i\in P$ we have an exceptional divisor $E_i\cong \p^{n-1}$, and
   $D'$ is disjoint from $E_i$. Thus  the
   pull-back of $D_P$ is linearly trivial on $E_i$, so $D_P\sim 0$.
     \qed

      \begin{thm}[Severi-Weil lemma]\label{weak.w-m.lem}
    Let $X$ be a geometrically normal, projective variety over an infinite field, $S$ a smooth variety and
 $D$ a Weil divisor on $S\times X$.
    Then there are finitely many  smooth, projective curves
   $\{C_i\subset X^{\rm sm}: i\in I\}$ such that, for all $s_1, s_2\in S$, 
$$
D_{s_1}\sim D_{s_2} \Leftrightarrow D_{s_1}|_{C_i}\sim D_{s_2}|_{C_i} \quad \forall i\in I.
$$
        \end{thm}

  Proof.
       Let $|H|$ be very ample on $X$.
  For any $C:=H_1\cap\cdots \cap H_{n-1}$ as in (\ref{weak.cl0.prop}.2)
  restriction to $C$ gives a morphism
  $\alpha_C: S\to \jac(C)$.
  Set $\Gamma_C:=S\times_{\jac(C)} S\subset S\times S$.

  By the Noetherian property, there is a finite set of curves
  $\{C_i: i\in I\}$ such that   $\cap_{i\in I}\Gamma_{C_i}$ equals the intersection of all the $\Gamma_C$.

  Thus, if $D_{s_1}|_{C_i}\sim  D_{s_2}|_{C_i}$ for every $i\in I$, then
  $D_{s_1}|_C\sim  D_{s_2}|_C$ for every $C=H_1\cap\cdots \cap H_{n-1}$, hence $D_{s_1}\sim  D_{s_2}$ by Lemma~\ref{weak.cl0.prop}.\qed

 \section{Further comments}

   I tried to simplify the proofs of Theorems~\ref{bert.irr.thm}--\ref{chev.thm}, thus several intermediate results   are stated in  weak forms, when this makes the proofs shorter.
   Here I point out some  more complete versions and  connections. These are not used in the paper, but should be useful in general.

 \begin{say}[Compactifications]\label{adv.com.say.1.}
  As a  refinement of (\ref{chev.thm.pf}.1), 
choose a finite,  affine cover $G=\cup_{i\in I} U_i$ and  embeddings
 $\tau_i:U_i\into \p^{n_i}$. The product gives a rational map
 $$
 G\map \times_{i\in I}U_i\into \times_{i\in I}\p^{n_i}.
 $$
 Let $X$ be the normalization of the closure of the image.
The advantage is that  $G\map  X$ is  an open embedding outside a
closed subset $Z\subset G$ of codimension $\geq 2$. Thus every divisor on $G$ is the restriction of a Weil divisor on $X$, and
$H^0(G, L)=H^0(G\setminus Z, L)$ for every line bundle $L$ on $G$. 

\smallskip

Note that the stronger compactification theorems do not quite work for our purposes.
  Nagata's   theorem \cite{MR0142549} says that any normal variety   has a normal, proper compactification; see     \cite[\href{https://stacks.math.columbia.edu/tag/0F41}{Tag 0F41}]{stacks-project}   for a modern proof.
  Weil's regularization theorem \cite{MR0074083} says that any rational $G$-action is birational to a regular $G$-action on a  proper, normal variety; see also \cite{MR337963, MR2562620}.
  However, these either do not guarantee projectivity, or appear to use
  Chevalley's theorem in order to prove  projectivity.

  We use the projectivity in Section~\ref{alb-pic.sect}  to find families of projective curves on $X$ that are contained in its smooth locus. There are proper, normal algebraic surfaces $S$ such that every curve  on $X$ passes through a singular point. 
\end{say}

 \begin{say}[Class group and Picard variety]\label{adv.com.say.2.}
The morphism $\alpha_{X, D}$ constructed in Theorem~\ref{weak.w-m.lem}
sends $S$ to an Abelian variety that depends on various choices.

Over arbitrary fields
            Matsusaka proved that there is a universal choice for the target,  called the   Picard variety of $X$; see \cite{MR62470, MR62471}.    I denote it by $\pics^{\rm wm}(X)$ for Weil and  Matsusaka.

            The definition of a Picard variety was changed  in \cite[VI]{fga}, and since then  {\it Picard variety} almost always means the Grothendieck version.
            The 2 versions agree for smooth varieties in characteristic 0, but  can have different dimensions in general.
            I plan to discuss this in detail later.
\end{say}

 \begin{say}[Irreducible sets of generators]\label{adv.com.say.3.}
    There are other examples where geometrically irreducible sets of generators behave better than arbitrary sets of generators, see for example
            \cite[Sec.3]{MR1714736} and \cite[Ques.14]{MR2894633}. 
   \end{say}

    \begin{ack}
      Partial  financial support    was provided  by the Simons Foundation   (grant number SFI-MPS-MOV-00006719-02).
    \end{ack}
 

\def\cprime{$'$} \def\cprime{$'$} \def\cprime{$'$} \def\cprime{$'$}
  \def\cprime{$'$} \def\dbar{\leavevmode\hbox to 0pt{\hskip.2ex
  \accent"16\hss}d} \def\cprime{$'$} \def\cprime{$'$}
  \def\polhk#1{\setbox0=\hbox{#1}{\ooalign{\hidewidth
  \lower1.5ex\hbox{`}\hidewidth\crcr\unhbox0}}} \def\cprime{$'$}
  \def\cprime{$'$} \def\cprime{$'$} \def\cprime{$'$}
  \def\polhk#1{\setbox0=\hbox{#1}{\ooalign{\hidewidth
  \lower1.5ex\hbox{`}\hidewidth\crcr\unhbox0}}} \def\cdprime{$''$}
  \def\cprime{$'$} \def\cprime{$'$} \def\cprime{$'$} \def\cprime{$'$}
\providecommand{\bysame}{\leavevmode\hbox to3em{\hrulefill}\thinspace}
\providecommand{\MR}{\relax\ifhmode\unskip\space\fi MR }
\providecommand{\MRhref}[2]{%
  \href{http://www.ams.org/mathscinet-getitem?mr=#1}{#2}
}
\providecommand{\href}[2]{#2}

\bigskip

Princeton University, Princeton NJ 08544-1000,

\email{kollar@math.princeton.edu}

\end{document}